\documentclass[12pt,reqno, oneside]{amsart}

\newtheorem{theorem}{Theorem}
\newtheorem{lemma}{Lemma}

\theoremstyle{definition}
\newtheorem*{remark}{Remark}
\newtheorem*{remarks}{Remarks}

\renewcommand{\t}[1]{\frac{t^#1}{#1!}}

\title[Three proofs of the Goulden-Litsyn-Shevelev Conjecture]
{Three proofs of the Goulden-Litsyn-Shevelev Conjecture on a Sequence Arising in Algebraic Geometry}

\author{Brian Drake$^1$, Ira M. Gessel$^{2\dagger}$ \and Guoce Xin$^3$}
\address{$^{1,2}$Department of Mathematics, Brandeis University, Waltham, MA, USA,
02454}
\address{$^3$Center for Combinatorics, LPMC, Nankai University, Tianjin, 300071, P. R. China}
\email{$^1$bdrake@brandeis.edu, $^2$gessel@brandeis.edu,
$^3$gxin@nankai.edu.cn}
\date{\today}
\thanks{$^\dagger$Partially supported by NSF Grant DMS-0200596}
\begin{document}

\begin{abstract}
We prove and generalize a conjecture of Goulden, Litsyn, and Shevelev that
certain Laurent polynomials related to the solution of a functional equation
have  only odd negative powers.
\end{abstract}

\maketitle
\thispagestyle{empty}

\section{Introduction} Consider the functional equation
\begin{equation}\label{e-original}
k (1+\chi) \log (1+\chi)=(k+1) \chi-t .\end{equation}
It is easily seen that \eqref{e-original} has a unique solution
$\chi=\chi(t,k)$
as a formal power series in $t$, and that the coefficient of $t^n$  in
$\chi(t,k)$ is a polynomial in $k$ of degree $n-1$:
\begin{equation}
\label{e-chi1}
\begin{split}
\chi(t,k) = t + k\t2 + (3k^2-k)\t3+(15k^3 - 10 k^2 +2k)\t4\qquad\\
  + (105k^4-105k^3+40k^2 -6k)\t5+\cdots
\end{split}
\end{equation}
Here the coefficients are sequence A075856 in the Online Encyclopedia of Integer Sequences \cite{OEIS}.

If we set $k=1$  then the
coefficients of $t^n\!/n!$ are the Betti numbers of the moduli space
of $n$-pointed stable curves of genus 0, as shown by Keel \cite{Keel}. These
coefficients give sequence A074059 in  \cite{OEIS}:
\begin{equation}
\chi(t, 1) = t + \frac{t^2}{2!} + 2 \frac{t^3}{3!} + 7
\frac{t^4}{4!} + 34 \frac{t^5}{5!} + \cdots
\end{equation}  For an arbitrary positive integer $k$, there is an interpretation of $\chi(t,k)$ in
terms of configuration spaces given by Manin \cite[p.~197]{Manin}.
See also Goulden, Litsyn, and Shevelev \cite{GLS} and the references given there. The series
expansion of $\chi$ has also been studied from a combinatorial
perspective.  See Dumont and Ramamonjisoa \cite[Proposition 6]{DR},
where $\chi$ is shown to count
functional digraphs by improper edges, and Zeng \cite[Corollary 11]{Zeng}, where $\chi$
is shown to count certain trees by improper edges.  Some
related bijections have been studied by Chen and Guo \cite{CG}.

Following Goulden, Litsyn, and Shevelev \cite{GLS}, let us
define $\mu_l(n)$ to be the coefficient of $k^{n-l}t^n/n!$ in
$\chi(t,k)$. The expansion \eqref{e-chi1} suggests that
$\mu_1(n)=1\cdot 3\cdots (2n-3)$, and thus that
\[\sum_{n=1}^\infty \mu_1(n)\t n = 1-\sqrt{1-2t}.\]
Goulden et al.~proved more generally that
for each $l$, $\sum_{n=0}^\infty \mu_l(n) t^n/n!$ is a Laurent polynomial $M_l(u)$ in $u=\sqrt{1-2t}$.

The first few polynomials $M_l(u)$, given in \cite{GLS},
are reproduced below:
\begin{align*}
M_1(u) &= 1-u,\\
M_2(u) &=-\frac16u^{-1}+\frac12-\frac12 u +\frac16 u^2,\\
M_3(u) &= \frac{1}{72} u^{-3}-\frac18 u^{-1}+\frac29-\frac18
u+\frac{1}{72}u^{3},\\
M_4(u) &=-\frac{1}{432}u^{-5}+\frac{1}{72} u^{-3} - \frac{1}{20}
u^{-1}+\frac{1}{18}-\frac{1}{72} u-\frac{1}{144} u^3+\frac{1}{270}
u^4, \\
M_5(u) &=\frac{5}{10368} u^{-7}-\frac{5}{1728}
u^{-5}+\frac{43}{5760}u^{-3}-\frac{59}{4320}
u^{-1}\\
&\qquad+\frac{1}{90}-\frac{5}{1152} u+\frac{1}{405} u^2+\frac{1}{576} u^3-\frac{1}{270}
u^4+\frac{23}{17280}
u^5.\\
&{}
\end{align*}
In this paper we give three proofs of the following conjecture of Goulden et al.~\cite[Conjecture 1]{GLS}:

\begin{theorem}\label{thm1}
The negative powers of $u$ in $M_l(u)$ are all odd.
\end{theorem}

Our first proof uses a change of variables to express the generating function in $y$
for the Laurent polynomials $M_n(u)$ in terms of the power series in $y$ (with coefficients that are
Laurent polynomials in $u$)
\begin{equation*}
S = e^{-y/2} y u \sqrt{1 + \frac{2y}{u^2} T},
\end{equation*}
where $T$ is a power series in $y$ only.
It is easily seen that any odd power of $S$ contains only odd (but possibly negative)
powers of $u$, while any even power of $S$ contains no negative powers of $u$. This change of variables
is related to an expansion of the Lambert $W$ function around its branch point at $-e^{-1}$
(though our proof uses only formal power series).

Our second proof uses properties of formal Laurent series, and shows that the
Goulden-Litsyn-Shevelev property holds much more generally for solutions of the functional equation
\begin{equation}
H = t+k\Lambda(H)
\end{equation}
where $\Lambda(z)$ is an arbitrary power series of the form
$\lambda_2 z^2+\lambda_3z^3+\cdots$ with $\lambda_2\ne0$; the Goulden-Litsyn-Shevelev conjecture corresponds to the case
$\Lambda(z) = (1+z)\log(1+z) -z$.

In our third proof, we consider the same generalization as in the second proof, but we apply Lagrange inversion to obtain an explicit formula from which it is clear that there are no negative even powers of $u$.

\section{First Proof of Theorem \ref{thm1}}\label{firstproof}

When $\chi$ is expanded as a power series in $t$, the coefficient of
$t^n$ is a polynomial in~$k$ of degree $n-1$.  To reverse the coefficients
of these polynomials, we introduce a new variable $y$  and set
$\xi = \xi(t, y) = \chi(yt,y^{-1})$,
so that $\sum_{n=0}^\infty \mu_l(n) t^n/n!$ is the coefficient of $y^l$ in $\xi$.
As a power series in $t$,
\begin{multline*}\label{xipowert}
\qquad \xi(t, y) = yt + y\frac{t^2}{2!} + (3y-y^2)\frac{t^3}{3!} +
(15y-10y^2+2y^3)\frac{t^4}{4!}\\
 + (105y-105y^2+40y^3-6y^4)\frac{t^5}{5!}+\cdots\qquad
\end{multline*}
and $\xi$ satisfies the functional equation
\begin{equation}\label{e-fun}
(1+\xi)\log(1+\xi) = (1+y)\xi-ty^2.
\end{equation}
As we shall see in the second proof, equation \eqref{e-fun} has two
solutions as power series in $t$ and $y$, but $\xi$ is the unique solution
in which the constant term in $t$ is zero.

We can expand $\xi$ as a power series in $y$, and as shown in \cite{GLS},
in this expansion the coefficient of $y^l$ is a Laurent polynomial in
$u=\sqrt{1-2t}$.

Replacing $t$ with $\tfrac 12 (1-u^2)$ in  \eqref{e-fun} gives
\begin{equation}\label{e-fun2}
(1+\xi)\log(1+\xi) = (1+y)\xi-\tfrac 12 (1-u^2) y^2.
\end{equation}

We will express $\xi$ in terms of an auxiliary power series $S = S(u, y)$ satisfying
\begin{equation}\label{Sdef}
1 + y + \tfrac{1}{2}(1-u^2)y^2 = e^y \left(1 - {S^2}\!/{2}\right)
\end{equation} in which the sign of $S$ is chosen so that $S = yu
+ \cdots$. Solving \eqref{Sdef} for $S$ gives
\begin{equation*}
S = e^{-y/2} y u \sqrt{1 + \frac{2y}{u^2} T},
\end{equation*}
where $T$ is the power series in $y$ given by
\begin{equation*}
T = \frac{e^y - 1 - y - y^2/2}{y^3}.
\end{equation*}
It is clear that $S$ is a power series in $y$ with coefficients that are Laurent polynomials in~$u$.

Our proof of Theorem \ref{thm1} relies on the following
lemma.

\begin{lemma}\label{Slem}
For every nonnegative integer $n$, $S(u, y)^n$ contains no even negative powers of~$u$.
\end{lemma}

\begin{proof}
If $n=2j$ is even, then
\[
S^{2j} = e^{-jy}y^{2j}u^{2j}\left(1 +
\frac{2y}{u^2}T\right)^{j} ,
\]
which is $u^{2j}$ times a polynomial in $u^{-1}$ of degree $2j$, so
no negative powers of $u$ appear.  If $n=2j+1$ is odd, then
\begin{align*}
S^{2j+1}&=e^{-jy-y/2}y^{2j+1}u^{2j+1}
  \left(1+\frac{2y}{u^2}T\right)^{j+1/2}\\
&= e^{-jy-y/2}y^{2j+1}u^{2j+1} \displaystyle\sum_{m \geq 0}
\binom{j+\frac{1}{2}}{m}\frac{2^m y^m}{u^{2m}}T^m.
\end{align*}
The sum does not terminate, so negative powers of $u$ appear.
However, we have an odd power of $u$ times a Laurent polynomial in
$u^2$, so only odd powers of $u$ occur.
\end{proof}

Now we want to show that $1 + \xi$ is $e^y$ times a power series
in $S$. Then Theorem \ref{thm1}  will follow from  Lemma \ref{Slem}.

Let $1 + \xi = e^y e^{-G}$ so that  $G = y - \log(1+\xi)  =yu+\cdots$. Making
these substitutions in equation \eqref{e-fun2} and simplifying, we get
\begin{equation*}
e^{-G}(1+G) = 1 - \frac{S^2}{2},
\end{equation*}
which yields
\begin{equation}
\label{e-SG}
S=\sqrt{2-2e^{-G}(1+G)}  =G - G^2/3 + \cdots
\end{equation}
Then \eqref{e-SG} can be inverted to express $G$ as a power series in $S$,
so $1 + \xi$ is $e^y$ times a power
series in $S$.  This completes the first proof of Theorem
\ref{thm1}.

\begin{remark}
What led us to these substitutions? We used Maple to solve equation (\ref{e-fun}), which gave
\begin{equation}\label{Maplesolution}
\xi = \exp\bigl(W(-(1+y+ty^2)e^{-1-y})+1+y\bigr)-1.
\end{equation}
Here $W$ is the  Lambert $W$ function, which satisfies
$W(z)e^{-W(z)}=z$.  (See \cite{W} for more information about this function.)
Equation \eqref{Maplesolution} suggests that if we set $1+\xi=e^{y+H}$ then $H$
may be simpler than $\xi$. (We take $G=-H$ to make the initial coefficient of $G$
positive.)
Since $W(z)$ has a branch point of order 2 at $z=-e^{-1}$,
$W(-e^{-1}+z)$ can be expanded as a power series in  $\sqrt{z}$.  Replacing
$z$ with $2S^2e^{-1}$ gives a power series with rational coefficients:
\begin{equation*}\label{Wsub}
W(-(1-2S^2)e^{-1}) = -1 + 2S - \frac{4}{3}S^2 + \frac{11}{9}S^3 -
\frac{172}{135}S^4 + \cdots .
\end{equation*}
To apply this expansion to $\xi$ we need to find $S$ with
\begin{equation*}
-(1-2S^2)e^{-1} = -(1+y+ty^2)e^{-1-y}.
\end{equation*}
This equation is equivalent to (\ref{Sdef}) with
$t = \frac{1}{2}(1-u^2)$.
\end{remark}

\section{Second Proof of Theorem
\ref{thm1}}\label{secondproof}

In this section we give a different proof of a generalization of Theorem
\ref{thm1}. Throughout this section,
$\Phi(z) = 1+\phi_1 z +\phi_2 z^2+\cdots$ is a power series in $z$
with constant term~1 in which the coefficients $\phi_1,\phi_2, \dots$  are indeterminates,
and $R$ is the ring of polynomials in the $\phi_i$, $u$, and $u^{-1}$, where $u$ is another
indeterminate.

\begin{lemma}
\label{L-Phi}
The equation
\begin{equation}
\label{e-Phi}
F^2\Phi(F) -2y F+(1-u^2)y^2=0
\end{equation}
has two power series solutions $F=f_1y + f_2y^2+\cdots$ with no constant term,
given by $F=Z(u,y)$ and $F=Z(-u,y)$, where
\begin{multline}
\label{e-Zcoeffs}
Z(u,y)=(1-u)y+\frac{(1-u)^3}{2u}\phi_1 y^2 +
\left(\frac{(1-u)^5(1+5u)}{8u^3}\phi_1^2  + \frac{(1-u)^4}{2u}\phi_2\right)y^3\\
+\left(\frac{(1-u)^7(1+7u+16u^2)}{16u^5}\phi_1^3
  +\frac{(1-u)^6(1+6u)}{4u^3}\phi_1\phi_2
  +\frac{(1-u)^5}{2u}\phi_3
\right)y^4  +\cdots.
\end{multline}
Moreover, the coefficient of $y^n$ in $Z(u,y)$ is a Laurent polynomial in $u$.
\end{lemma}

\begin{proof}
Let $F=f_1y+f_2y^2+\cdots$. Substituting in \eqref{e-Phi} and equating coefficients of $y^2$ gives
$(f_1-1)^2-u^2=0$, so if $F$ satisfies \eqref{e-Phi} then $f_1$ is either $1-u$ or $1+u$.
The coefficient of $y^{n+1}$ in $F^2\Phi(F)$ is a polynomial in  $\phi_1,\dots, \phi_{n-1},f_1,\dots, f_n$
in which $f_n$ occurs only in the term $2f_1f_n$,
and the coefficient of $y^{n+1}$ in $-2y F$ is $-2f_{n-1}$.
Thus for $n\ge2$, the coefficient of $y^{n+1}$ in the left side of \eqref{e-Phi} is a
polynomial in $\phi_1,\dots,
\phi_{n-1},f_1,\dots, f_n$ in which
$f_n$ occurs only in the term $2(f_1-1)f_n$.
Thus after a choice of either $f_1= 1-u$ or $f_1=1+u$,
there is a unique solution of \eqref{e-Phi} in which $f_n$ is a polynomial in $\phi_1,\dots,
\phi_{n-1},f_1,\dots, f_{n-1}$ divided by
$f_1-1$. Since $f_1 -1$ is either $u$ or $-u$, it follows that
$f_n$ is a Laurent polynomial in $u$. Replacing $u$ by $-u$ switches the two possibilities for $f_1$, so it must also
switch the two solutions.
\end{proof}

In the next lemma we will work in the ring $R((X))[[y]]$ of formal power series in $y$ with coefficients that
are Laurent series in $X$. These are series of the form
\[\sum_{i=-\infty}^\infty\sum_{j=0}^\infty a_{ij}X^i y^j\]
such that for each $j\ge0$, $a_{ij}=0$ for all but finitely many negative values of $i$.
Related applications of formal Laurent series can be found in \cite{G} and \cite{X}.

\begin{lemma}
 Let $X$ be an indeterminate and let $Z(u,y)$ be as in Lemma \ref{L-Phi}. Then in
the formal Laurent series ring $R((X))[[y]]$, we have
\begin{equation}
\label{e-factor}
\Phi(X) -\frac{2y}{X} +\frac{(1-u^2)y^2}{X^2}
=\left(1-\frac{Z(u,y)}{X}\right) \left(1-\frac{Z(-u,y)}{X}\right) T,
\end{equation}
where $T\in R[[X,y]]$ and $T$ has constant term 1 in $X$ and $y$.
\end{lemma}

\begin{proof}
Let
\[X^2\Phi(X) -2yX+(1-u^2)y^2 =\sum_{j=0}^\infty c_jX^j\]
and let
\begin{equation*}
T_1=\left(1-\frac{Z(u,y)}{X}\right)^{-1} \left(\Phi(X) -\frac{2y}{X} +\frac{(1-u^2)y^2}{X^2}\right).
\end{equation*}
In the ring  $R((X))[[y]]$ we have the expansion
\begin{equation*}
\left(1-\frac{Z(u,y)}{X}\right)^{-1} = \sum_{n=0}^\infty \left(\frac{Z(u,y)}{X}\right)^{n}.
\end{equation*}
Then writing $Z$ for $Z(u,y)$, we have
\[
T_1=\sum_{i=0}^\infty
\left(\frac ZX\right)^{\!\! i}\,\sum_{j=0}^\infty c_j X^{j-2},
\] so the coefficient of
$X^{-m}$ in $T_1$ is
\[\sum_{j\ge\max(2-m,0)} c_jZ^{m+j-2}.\]

Thus for $m\ge2$ the coefficient of $X^{-m}$ in $T_1$
is  0.
A similar argument applied to $T=\left(1-Z(-u,y)/X\right)^{-1}T_1$ shows that the coefficient of $X^{-m}$
in $T$ is 0 for $m\ge1$, so $T\in R[[X,y]]$.

Finally, setting $y=0$ in \eqref{e-factor} shows that the constant term in $y$ in $T$ is
$\Phi(X)$, so the constant term in $X$ and $y$ is 1.

\end{proof}

\begin{lemma}
\label{L-noeven}
For any nonnegative integer $r$, $Z(u,y)^r$ contains no even negative powers of~$u$.
\end{lemma}
\begin{proof}
It is sufficient to prove  that $Z(u,y)^r+Z(-u,y)^r$ contains no negative powers of~$u$. Let $\Psi$ be the expression appearing in  \eqref{e-factor}. Since the left side of \eqref{e-factor} has constant term 1 and has no negative powers of
$u$, $\log \Psi$ is a well-defined element of
$R((X))[[y]]$ with no negative powers of $u$.

The logarithm of the right side of \eqref{e-factor} is
\begin{equation*}
-\sum_{r=1}^\infty
\frac{Z(u,y)^r+Z(-u,y)^r}{rX^r} +\log T
\end{equation*}
Since $T$ contains no negative powers of $X$, for $r>0$ the coefficient of $X^{-r}$ in
$\log\Psi$ is
$-\tfrac{1}{r} \left(Z(u,y)^r+Z(-u,y)^r\right)$, which therefore
contains no negative powers of $u$.
\end{proof}

We can now prove our generalization of the
Goulden-Litsyn-Shevelev Conjecture.

\begin{theorem}
\label{thm2}
Let $\Lambda(z)=\lambda_2 z^2+\lambda_3 z^3+\cdots$ be a formal power series, where $\lambda_2\ne0$.  Then the equation
\begin{equation}\label{e-H}
H = t+k\Lambda(H)
\end{equation}
has a unique power series solution $H(t,k) =
t+\lambda_2 kt^2 +(\lambda_3 k +2\lambda_2^2 k^2)t^3+\cdots$,
in which the coefficient of $t^n$ is a polynomial in $k$ of degree $n-1$.
Let $\mu_l(n)/n!$ be the coefficient of $k^{n-l}t^n$ in $H(t,k)$.
Then $\sum_{n=0}^\infty \mu_l(n)t^n/n!$ is a Laurent polynomial in $\sqrt{1-4\lambda_2 t}$
with no even negative powers.\footnote{The $n!$ in the denominator is only for compatibility
with the original Goulden-Litsyn-Shevelev conjecture.}
\end{theorem}

\begin{proof} We shall prove only the case $\lambda_2=\tfrac 12$ of the theorem. The general case then follows,
after a short calculation, by
replacing $k$ with $2\lambda_2 k$ in \eqref{e-H}.

We prove only the last assertion of the theorem; the other assertions are straightforward.
Set $J(t,y)=H(yt, y^{-1})$ so that
$\sum_{n=0}^\infty \mu_l(n)t^n/n!$ is the coefficient of $y^l$ in $J$ and set
$\Phi(z)=2\Lambda(z)/z^2 = 1+2\lambda_3z+\dots$. Then \eqref{e-H} becomes
$J = yt+\tfrac12 y^{-1} J^2\Phi(J)$,
which may be rewritten as
\begin{equation}
\label{e-J}
J^2\Phi(J) -2yJ+2ty^2 =0.
\end{equation}
By Lemma \ref{L-Phi}, $J=Z(\sqrt{1-2t},y)$ is a solution of \eqref{e-J} with no constant term in
$t$, so it must be the unique such solution, and the conclusion follows from the case $m=1$ of Lemma
\ref{L-noeven}.

\end{proof}

\begin{remark}
The reader may wonder how the two solutions of Lemma \ref{L-Phi} became the unique solution of
Theorem \ref{thm2}. We have two solutions, $Z(\sqrt{1-2t},y)$ and  $Z(-\sqrt{1-2t},y)$ of \eqref{e-J}.
The coefficient of  $t^n$ in  $Z(\sqrt{1-2t},y)$ is a polynomial in $y$, but the coefficient of $t^n$ in
$Z(-\sqrt{1-2t},y)$ is not a polynomial in $y$, so although $Z(-\sqrt{1-2t},y)$ corresponds to a solution
(with nonzero constant term in $t$) of \eqref{e-J}, it does not correspond to a formal power series solution of
\eqref{e-H}.
\end{remark}

Theorem \ref{thm1} is obtained from  Theorem \ref{thm2} by taking
$\Lambda(z) = (1+z)\log(1+z) -z$.

Two cases of Theorem \ref{thm2} with  simple combinatorial interpretations are worth mentioning. First take
$\Lambda(z)= e^z -z -1$. (It is interesting to note that in this case, as in \eqref{e-original},
the solution can be expressed in terms of the Lambert $W$ function.)Here
\begin{multline*}
\qquad H(t,k) =t+k\frac{t^2}{2!}+(k+3k^2)\frac{t^3}{3!}+(k+10k^2+15k^3)\frac{t^4}{4!}\\
  +(k+25k^2+105k^3+105k^4)\frac{t^5}{5!}+\cdots\qquad
\end{multline*}
The coefficient of $k^it^n/n!$ in $H(t,k)$ is the number of rooted trees with $i$
unlabeled internal vertices, each with at least two children, and $n$ labeled leaves. (These are sometimes called \emph{phylogenetic
trees}.) Setting $k=1$ gives
\[
 H(t,1) = t+\frac{t^2}{2!}+4\frac{t^3}{3!}+26\frac{t^4}{4!}+236\frac{t^5}{5!}+
  2752\frac{t^6}{6!}+39208\frac{t^7}{7!}+660032\frac{t^8}{8!}\cdots.
\]
These coefficients are sequence A000311 in \cite{OEIS}.

Another example is $\Lambda(z)= {z^2}/{(1-z)}$, so
\begin{equation}
\label{e-trees}
H = t +k \frac{H^2}{1-H}.
\end{equation}
Here we can solve for $H$ explictly, obtaining
\begin{equation}
\label{e-Hex}
H=\frac{1+t -\sqrt{(1-t)^2 -4kt}}{1+k} = t+\sum_{n=2}^{\infty}\sum_{i=1}^{n-1}
\frac 1i \binom{n-2}{ i-1}\binom{n+i-1}{i-1}k^i t^n,
\end{equation}
where the formula for the coefficients is easily derived by Lagrange inversion.
It is clear from \eqref{e-trees} that the coefficient of $k^i t^n$ in $H(t,k)$ is the number of
ordered trees with $n$~leaves and
$i$ internal vertices, in which every internal vertex has at least two children. These coefficients are
sequence A033282 in
\cite{OEIS}, which gives many references and some additional combinatorial interpretations, notably
in terms of dissections of a polygon.

Using the explicit formula \eqref{e-Hex} for the solution of \eqref{e-trees}  we can
verify this case of Theorem
\ref{thm2} directly. In fact, a stronger statement holds here: if $\rho_l(n)$ is the coefficient of
$k^{n-l} t^n$ in $H(t,k)$ (here it is appropriate to use ordinary, rather than exponential, generating
functions) then
$\sum_{n=0}^\infty \rho_l(n) t^n$ is a Laurent polynomial in $\sqrt{1-4t}$ in which the only even powers of
$\sqrt{1-4t}$ are $(\sqrt{1-4t})^0$ and $(\sqrt{1-4t})^2$. To see this, we note that it follows from \eqref{e-Hex}
that
\begin{equation*}
\sum_{l=0}^\infty y^l \sum_{n=0}^\infty \rho_l(n) t^n =
 H(yt, y^{-1})= y\frac{1+yt-\sqrt{(1-yt)^2-4t}}{2(1+y)}=J(t,y)
\end{equation*}
Setting $u=\sqrt{1-4t}$, so that $t=(1-u^2)/4$, we have
\begin{equation*}
J((1-u^2)/4, y)=y\frac{4+y(1-u^2) -u\sqrt{16+8y(1-u^{-2})+y^2(u-u^{-1})^2}}{8(1+y)}.
\end{equation*}
Note that the expression under the square root sign involves only even powers of $u$, and it is multiplied by $u$,
giving only odd powers of $u$. Thus the contribution to even powers of $u$ is
\begin{equation*}
y\frac{4+y(1-u^2)}{8(1+y)} = \frac y2 -\frac{y^2(3+u^2)}{8(1+y)} = \frac y2 + \sum_{l=2}^\infty \tfrac18 (-1)^{l-1}
(3+u^2) y^l.
\end{equation*}

\section{Third proof of Theorem \ref{thm1}}
\label{thirdproof}

In our third proof, we use Lagrange inversion to prove Lemma \ref{L-noeven} (from which, as we have seen, Theorem~\ref{thm2} follows easily) by giving an explicit formula for $Z(u,y)^r$ that makes it clear that it has no even negative powers of $u$.

\begin{theorem}
\label{T-exp}
With $Z(u,y)$ as in Lemma \ref{L-Phi}, we have for any positive integer $r$,
\begin{multline}
\label{e-Z}
Z(u,y)^r = y^r(1-u)^r
+ \sum_{n=r+1}^\infty y^n \\
\times\!\!\sum_{m_1+2m_2+\cdots = n-r}
(-1)^m\frac{r}{2m+n}\binom{m}{m_1,m_2,\cdots}P_{m,n}(u)
\phi_1^{m_1}\phi_2^{m_2}\cdots,\qquad
\end{multline}
where $m=m_1+m_2+\cdots$ and
\begin{equation}
\label{e-explicit}
P_{m,n}(u)  = \sum_{i=1-2m}^n \binom{2m+n}{2m+i}
\binom{m+i/2}{m}(-u)^i.
\end{equation}
\end{theorem}

\begin{proof}
Let us set $\Psi(z) = z^2(\phi_1 z+\phi_2 z^2+\cdots),$ so that with
$\Phi(z)$ as defined at the beginning of Section
\ref{secondproof} we have $z^2\Phi(z)=z^2+\Psi(z)$. Then \eqref{e-Phi} may be written
$(F-y)^2 = u^2 y^2 - \Psi(F)$. Taking square
roots gives
\begin{equation}
\label{e-Lag1}
F-y = -uy\sqrt{
\vrule height 0pt width 0pt depth8pt
\smash[b]{1-\dfrac{\Psi(F)}{u^2y^2}}},
\end{equation}
where the sign of the square root is chosen so as to give the solution in which the coefficient of $y$ is
$1-u$ rather than $1+u$. To apply  Lagrange inversion, we must make a slight modification in  \eqref{e-Lag1}.
It is not difficult to show, e.g., by equating coefficients in~\eqref{e-Phi}, that the coefficient of $y^n$ in
$Z(u,y)$ is a sum of terms in $\phi_1^{m_1}\phi_2^{m_2}\cdots
\phi_{n-1}^{m_{n-1}}$, where $m_1+2 m_2 +\cdots (n-1) m_{n-1}=n-1$. Thus the variable $y$ in $Z(u,y)$ is
redundant; $Z(u,y)$ can be recovered from $Z(u,1)$ be replacing each $\phi_i$ with $y^i\phi_i$ and then
multiplying by $y$. So it is enough to solve $F=1-u\sqrt{1-u^{-2}\Psi(F)}$ for $F$ as a power series in $\phi_1$,
$\phi_2$, \dots. Finally, to put this equation into a form to which Lagrange inversion, in its usual form, can be
applied (see, e.g., Stanley \cite[Theorem 5.4.2]{EC2}) we introduce a new redundant variable $x$ and consider the
equation
\begin{equation}
\label{e-Lag2}
F=x\bigl(1-u\sqrt{1-u^{-2}\Psi(F)}\,\bigr).
\end{equation}
for $F$ in as power series in 
$x$.
Applying Lagrange inversion, we have for any positive integer
$r$,
\begin{align*}
[x^s] F^r&=\frac rs [t^{s-r}]\bigl(1-u\sqrt{1-u^{-2}\Psi(t)}\bigr)^s\\
  &= \frac rs [t^{s-r}]\sum_{j=0}^s \binom{s}{j}(-u)^j \bigl(1-u^{-2}
        \Psi(t)\bigr)^{j/2}\\
  &=\frac rs [t^{s-r}]\sum_{j=0}^s \binom{s}{j}(-u)^j
       \sum_m \binom{j/2}m (-u^{-2})^m\Psi(t)^m\displaybreak[0]\\
  &=\frac rs [t^{s-r}]\sum_{j=0}^s \binom{s}{j}\sum_m (-1)^m
      \binom{j/2}{m} (-u)^{j-2m}\\
   &\qquad\quad \times\!\!\sum_{m_1+m_2+\cdots = m}\binom{m}{m_1, m_2,\cdots}
      \phi_1^{m_1}\phi_2^{m_2}\cdots t^{2m+m_1+
      2m_2+\cdots}.\\
\end{align*}
Thus with $m=m_1+m_2+\cdots$ and $n=r+m_1+2m_2+\cdots$, setting $s=2m+n$ gives
\begin{multline*}
\qquad\qquad\qquad F^r = \sum_{m_1, m_2,\dots}x^{2m+n}\frac{r}{2m+n}
      \binom{m}{m_1,m_2,\cdots}\phi_1^{m_1}\phi_2^{m_2}\cdots \\
   \times(-1)^m\sum_{j=0}^{2m+n}\binom{2m+n}j\binom{j/2}m (-u)^{j-2m}\qquad\qquad\qquad
\end{multline*}
To obtain \eqref{e-Z}, we first set $x=1$, replace each $\phi_i$ with $y^i\phi_i$, and multiply by $y^r$.
The contribution from $m_1=m_2=\cdots = 0$ (which gives $m=0$ and $n=r$) is
\[y^r\sum_{j=0}^r\binom rj (-u)^j = y^r(1-u)^r.\]
If $m>0$ then  $\binom {j/2}m=0$ for $j=0$, so we may start the inner sum on $j$ with $j=1$, rather than $j=0$. Finally, setting $j=i+2m$ yields \eqref{e-Z}.
\end{proof}

It follows immediately from Theorem \ref{T-exp} that $Z(u,y)$ has no negative even powers of~$u$, since $\binom{m+i/2}{m}$ is 0 for $i$ even and negative.

\begin{remarks}
The coefficients of $Z(u,y)$ as displayed in \eqref{e-Zcoeffs} show divisibility  by powers of $1-u$, and it is not difficult to prove from \eqref{e-explicit} that the numerator of $P_{m,n}(u)$ is divisible by $(1-u)^{m+n}$.

Theorem \ref{T-exp} can be generalized to the equation $(F-y)^p + \Psi(F) - u^p y^p=0$.
The coefficients are Laurent polynomials in $u$ in which the coefficient of $u^{-ip}$ is 0 for every positive integer $i$.
\end{remarks}

\end{document}